# Эволюционные выводы энтропийной модели расчета матрицы корреспонденций

*Гасников А.В., Гасникова Е.В., Мендель М.А., Чепурченко К.В.*

Лаборатория структурных методов анализа данных в предсказательном моделировании (ПреМоЛаб) Факультета управления и прикладной математики Национального исследовательского Университета «Московский физико-технический институт»
Институт проблем передачи информации им. А.А. Харкевича Российской академии наук

В работе приводятся два эволюционных способа вывода классической энтропийной модели расчета матрицы корреспонденций. Оба способа базируются на подходе стохастической химической кинетики в условиях детального баланса. Первый способ базируется на описании бинарных реакций обмена местами жительства жителей города, второй способ базируется на унарных реакциях, составляющих фундамент стохастической части современной теории популяционных игр. Оба способа позволяют глубже понять физический смысл ряда параметров энтропийной модели, и получить ответы на несколько открытых вопросов из смежных областей. Например, приведенные в статье результаты позволяют обосновать способ отбора единственного равновесия в модели равновесного распределения потоков по путям в случае, когда равновесий много.

## 1. Введение

Настоящая статья предлагает два способа объяснения популярной в урбанистике энтропийной модели расчета матрицы корреспонденций [1]. Обе предложенные модели в своей основе имеют марковский процесс, живущий в пространстве огромной размерности (говорят, что такой процесс порождает макросистему). Более точно, этот марковский процесс представляет собой ветвящийся процесс специального вида “модель стохастической химической кинетики” [2, 3]. Более того, в обоих случаях имеет место условие детального равновесия [3–5]. Из общих результатов [3–5] отсюда сразу можно заключить, что инвариантная (стационарная) мера такого процесса будет иметь вид мультиномиального распределения, сосредоточенного на аффинном многообразии, определяемым линейными законами сохранения введенной динамики. Исходя из явления концентрация меры [6] можно ожидать концентрацию этой меры около наиболее вероятного состояния с ростом размера макросистемы. Такое состояние, естественно, и принимать за равновесие изучаемой макросистемы, поскольку с большой вероятностью на больших временах мы найдем систему в малой окрестности такого состояния. Поиск равновесия сводится по теореме Санова [7] к задаче энтропийно-линейного программирования (3) [8]. Собственно, именно таким образом и планируется объяснить, почему для описания равновесной матрицы корреспонденций необходимо решить задачу ЭЛП.

В п. 2 на базе бинарных реакций обменного типа, популярных в различного рода физических и социально-экономических приложениях моделей стохастической химической кинетики [9, 10], будет приведен первый способ эволюционного вывода энтропийной модели расчета матрицы корреспонденций. Второй способ будет описан в п. 3. Он базируется на классических понятиях популяционной теории игр [11]. Например, таких как логит динамика и игра загрузки (именно к такой игре сводится поиск равновесного распределения потоков по путям). Второй способ также базируется на редукции задачи расчета матрицы корреспонденций к задаче поиска равновесного распределения потоков по путям. Нетривиальным и, по-видимому, новым здесь является эволюционно-экономическая интерпретация двойственных множителей, обобщающая известную стационарную конструкцию, восходящую к Л.В. Канторовичу.

Отметим, что результаты данной работы являются обобщением результатов работ [12, 13]. Обобщение заключается в том, что в формулировках основных результатов (теоремы 1 и 2) фигурируют точные (не улучшаемые) оценки времени сходимости изучаемых макросистем на равновесие и плотности концентрации инвариантной меры в окрестности равновесий.

Мы намеренно опускаем вопросы численного решения, возникающих задач ЭЛП. Подробнее об этом можно посмотреть, например [8, 14].

## 2. Вывод на основе обменов

Приведем, базируясь на работе [12], эволюционное обоснование одного из самых популярных способов расчета матрицы корреспонденций, имеющего более чем сорокалетнюю историю, – энтропийной модели [1].

Пусть в некотором городе имеется $n$ районов, $L_i > 0$ – число жителей $i$-го района, $W_j > 0$ – число работающих в $j$-м районе. При этом $N = \sum_{i=1}^{n} L_i = \sum_{j=1}^{n} W_j$, – общее число жителей города, $n^2 \ll N$. Далее под $L_i \geq 0$ будет пониматься число жителей района, выезжающих в типичный день за рассматриваемый промежуток времени из $i$-го района, а под $W_j \geq 0$ – число жителей города, приезжающих на работу в $j$-й район в типичный день за рассматриваемый промежуток времени. Обычно, так введенные, $L_i$, $W_j$ рассчитываются через число жителей $i$-го района и число работающих в $j$-м районе с помощью более менее универсальных (в частности, не зависящих от $i$, $j$) коэффициентов пропорциональности. Эти величины являются входными параметрами модели, т.е. они не моделируются (во всяком случае, в рамках выбранного подхода). Для долгосрочных расчетов с разрабатываемой моделью требуется иметь прогноз изменения значений этих величин.

Обозначим через $d_{ij}(t) \geq 0$ – число жителей, живущих в $i$-м районе и работающих в $j$-м в момент времени $t$. Со временем жители могут только меняться квартирами, поэтому во все моменты времени $t \geq 0$

$$d_{ij}(t) \geq 0,\ \sum_{j=1}^{n} d_{ij}(t) \equiv L_i,\ \sum_{i=1}^{n} d_{ij}(t) \equiv W_j,\ i, j = 1, ..., n.$$

Определим

$$\mathrm{A} = \left\{ d_{ij} \geq 0:\ \sum_{j=1}^{n} d_{ij} = L_i, \sum_{i=1}^{n} d_{ij} = W_j, i, j = 1, ..., n \right\}.$$

Опишем основной стимул к обмену: работать далеко от дома плохо из-за транспортных издержек. Будем считать, что эффективной функцией затрат [12] будет $R(T) = \beta T/2$, где $T > 0$ – время в пути от дома до работы (в общем случае под $T$ стоит понимать затраты, в которые может входить не только время), а $\beta > 0$ – настраиваемый параметр модели (который также можно проинтерпретировать и даже оценить, что и будет сделано ниже).

Теперь опишем саму динамику. Пусть в момент времени $t \geq 0$ $r$-й житель живет в $k$-м районе и работает в $m$-м, а $s$-й житель живет в $p$-м районе и работает в $q$-м. Тогда $\lambda_{k,m;\,p,q}(t)\Delta t + o(\Delta t)$ – есть вероятность того, что жители с номерами $r$ и $s$ ( $1 \leq r < s \leq N$ ) "поменяются" квартирами в промежутке времени $(t, t+\Delta t)$. Вероятность обмена местами жительства зависит только от мест проживания и работы обменивающихся:

$$\lambda_{k,m;\,p,q}(t) \equiv \lambda_{k,m;\,p,q} = \lambda N^{-1} \exp\Big( \underbrace{R(T_{km}) + R(T_{pq})}_{\text{суммарные затраты до обмена}} - \underbrace{\left(R(T_{pm}) + R(T_{kq})\right)}_{\text{суммарные затраты после обмена}} \Big) > 0,$$

где коэффициент $0 < \lambda = \mathrm{O}(1)$ характеризует интенсивность обменов. Совершенно аналогичным образом можно было рассматривать случай "обмена местами работы". Здесь стоит оговориться, что "обмены" не стоит понимать буквально – это лишь одна из возможных интерпретаций. Фактически используется, так называемое, "приближение среднего поля” [2, 5], т.е. некое равноправие агентов (жителей) внутри фиксированной корреспонденции и их независимость.

Согласно эргодической теореме для марковских цепей (в независимости от начальной конфигурации $\{d_{ij}(0)\}_{i=1,\, j=1}^{n,\, n}$) [3, 9–12, 15, 16] предельное распределение совпадает со стационарным (инвариантным), которое можно посчитать (получается проекция прямого произведение распределений Пуассона на А):

$$\lim_{t\to\infty} P\left(d_{ij}(t)=d_{ij},\, i,\, j=1,\ldots,n\right) = Z^{-1}\prod_{i,j=1}^{n}\exp\left(-2R\left(T_{ij}\right)d_{ij}\right)\cdot\left(d_{ij}!\right)^{-1} \overset{def}{=} p\left(\{d_{ij}\}_{i=1,\, j=1}^{n,\, n}\right), \quad (1)$$

где $\{d_{ij}\}_{i=1,\, j=1}^{n,\, n} \in \mathrm{A}$, а “статсумма” $Z$ находится из условия нормировки получившейся “пуассоновской” вероятностной меры. Отметим, что стационарное распределение $p\left(\{d_{ij}\}_{i=1,\, j=1}^{n,\, n}\right)$ удовлетворяет условию детального равновесия [5, 11]:

$$\left(d_{km}+1\right)\left(d_{pq}+1\right)p\left(\{d_{11},\ldots,d_{km}+1,\ldots,d_{pq}+1,\ldots,d_{pm}-1,\ldots,d_{kq}-1,\ldots,d_{nn}\}\right)\lambda_{k,m;\,p,q} =$$

$$= d_{pm}d_{kq}\,p\left(\{d_{ij}\}_{i=1,\, j=1}^{n,\, n}\right)\lambda_{p,m;\,k,q}.$$

При $N \gg 1$ распределение $p\left(\{d_{ij}\}_{i=1,\, j=1}^{n,\, n}\right)$ экспоненциально сконцентрировано на множестве А в $\mathrm{O}\left(\sqrt{N}\right)$ окрестности наиболее вероятного значения $d^* = \{d_{ij}^*\}_{i=1,\, j=1}^{n,\, n}$, которое определяется, как решение задачи энтропийно-линейного программирования [3, 12, 13]:

$$\ln p\left(\{d_{ij}\}_{i=1,\, j=1}^{n,\, n}\right) \sim -\sum_{i,j=1}^{n} d_{ij}\ln\left(d_{ij}\right) - \beta\sum_{i,j=1}^{n} d_{ij}T_{ij} \to \max_{\{d_{ij}\}_{i=1,\, j=1}^{n,\, n}\in(\mathrm{A})}. \quad (2)$$

Это следует из теоремы Санова о больших уклонениях для мультиномиального распределения [7] (на распределение (1) можно смотреть также как на проекцию мультиномиального распределения на А)

$$\frac{N!}{d_{11}!\cdot\ldots\cdot d_{ij}!\cdot\ldots\cdot d_{nn}!}\, p_{11}^{d_{11}}\cdot\ldots\cdot p_{ij}^{d_{ij}}\cdot\ldots\cdot p_{nn}^{d_{nn}} = \exp\left(-N\sum_{i,j=1}^{n}\nu_{ij}\ln\left(\nu_{ij}/p_{ij}\right)+\bar{R}\right),$$

где $\nu_{ij} = d_{ij}/N$, $|\bar{R}| \le \frac{n^2}{2}\left(\ln N + 1\right)$, и формулы Тейлора (с остаточным членом второго порядка в форме Лагранжа), примененной к энтропийному функционалу в точке $d^*$, заданному на А [1].

Сформулируем более точно полученный результат.

**Теорема 1.** *Для любого* $d(0)\in \mathrm{A}$ *существует такая константа* $c_n\left(d(0)\right) > 0$, *что для всех* $\sigma \in (0, 0.5)$, $t \ge c_n\left(d(0)\right)\ln N$ *имеет место неравенство*

$$P\left(\frac{\left\|d(t)-d^*\right\|_2}{N} \ge \frac{2\sqrt{2}+4\sqrt{\ln\left(\sigma^{-1}\right)}}{\sqrt{N}}\right) \le \sigma,$$

*где* $d(t)=\left\{d_{ij}(t)\right\}_{i=1,\ j=1}^{n,\ n} \in \mathrm{A}$ .

**Схема доказательства.** Установим сначала оценку для плотности концентрации стационарной меры

$$\lim_{t\to\infty} P\left(d_{ij}(t)=d_{ij},\, i,\, j=1,\ldots,n\right)=\frac{N!}{d_{11}!\cdot\ldots\cdot d_{ij}!\cdot\ldots\cdot d_{nn}!}\, p_{11}^{d_{11}}\cdot\ldots\cdot p_{ij}^{d_{ij}}\cdot\ldots\cdot p_{nn}^{d_{nn}}\, .$$

Из неравенства Хефдинга в гильбертовом пространстве [6] следует ( $\varepsilon \ge \sqrt{2/N}$ )

$$\lim_{t\to\infty} P\left(\left\|d(t)-d^*\right\|_2 \ge \varepsilon N\right) \le \exp\left(-\frac{1}{4N}\left(\varepsilon N-\sqrt{2N}\right)^2\right).$$

Беря в этом неравенстве

$$\varepsilon=\frac{\sqrt{2}+2\sqrt{\ln\left(\sigma^{-1}\right)}}{\sqrt{N}}\,,$$

получим

$$\lim_{t\to\infty} P\left(\frac{\left\|d(t)-d^*\right\|_2}{N} \ge \varepsilon\right) \le \sigma.$$

Однако если не переходить к пределу по времени, а лишь обеспечить при $t \ge T_n\left(\varepsilon, N; d(0)\right)$ выполнение условия

$$E\left\|\frac{d(t)-d^*}{N}\right\|_2 \le \frac{\varepsilon}{2}\,,$$

то при $t \ge T_n\left(\varepsilon, N; d(0)\right)$

$$P\left(\left\|d(t)-d^*\right\|_2 \ge \varepsilon N\right) \le \exp\left(-\frac{1}{4N}\left(\varepsilon N-\left(\sqrt{2N}+\frac{\varepsilon N}{2}\right)\right)^2\right)=\exp\left(-\frac{1}{4N}\left(\frac{\varepsilon N}{2}-\sqrt{2N}\right)^2\right).$$

Беря в этом неравенстве

$$\varepsilon=\frac{2\sqrt{2}+4\sqrt{\ln\left(\sigma^{-1}\right)}}{\sqrt{N}}\,,$$

получим при $t \ge T_n\left(\varepsilon, N; d(0)\right)$ аналогичное неравенство

$$P\left(\frac{\left\|d(t)-d^*\right\|_2}{N} \ge \varepsilon\right) \le \sigma.$$

Отметим, что с точностью до мультипликативной константы эта оценка не может быть улучшена. Это следует из неравенства Чебышёва

$$P\left(X \geq EX - \varepsilon\right) \geq 1 - \frac{\operatorname{Var}\left(X\right)}{\varepsilon^2},$$

при

$$X = \frac{\left\| d\left(t\right) - d^* \right\|_2}{N},\ \ d^* = N \cdot \left(n^{-2}, \ldots, n^{-2}\right)^T,\ \ N \gg n^2.$$

Осталось оценить зависимость $T_n\left(\varepsilon, N; d\left(0\right)\right)$. Для этого воспользуемся неравенством Чигера [16]. Поставим в соответствие нашему марковскому процессу его дискретный аналог (в дискретном времени с шагом $\lambda N^{-1}$) – случайное блуждание (со скачками, соответствующими парным реакциям, введенным выше) на целочисленных точках части гиперплоскости, задаваемой множеством А. Граф, на котором происходит блуждание, будем обозначать $G = \left\langle V_G, E_G \right\rangle$. Пусть $\pi\left(\cdot\right)$ стационарная мера этого блуждания (в нашем случае – мультиномиальная мера, экспоненциально сконцентрированая на множестве А в $\mathrm{O}\left(\sqrt{N}\right)$ окрестности наиболее вероятного значения $d^*$), а $P = \left\| p_{ij} \right\|_{i,j=1}^{|V_G|,|V_G|}$ – матрица переходных вероятностей. Тогда, вводя константу Чигера ($\bar{S} = V_G \setminus S$)

$$h\left(G\right) = \min_{S \subseteq V_G : \pi(S) \leq 1/2} P\left(S \to \bar{S} \middle| S\right) = \min_{S \subseteq V_G : \pi(S) \leq 1/2} \frac{\sum\limits_{(i,j) \in E_G : i \in S,\ j \in \bar{S}} \pi\left(i\right) p_{ij}}{\sum\limits_{i \in S} \pi\left(i\right)}$$

и время выхода блуждания (стартовавшего из $i \in V_G$) на стационарное распределение

$$T\left(i, \varepsilon\right) = \mathrm{O}\left(\lambda N^{-1} \cdot h\left(G\right)^{-2} \left(\ln\left(\pi\left(i\right)^{-1}\right) + \ln\left(\varepsilon^{-1}\right)\right)\right),$$

получим для любых $i = 1, \ldots, |V_G|$, $t \geq T\left(i, \varepsilon\right)$

$$\left\| P^t\left(i, \cdot\right) - \pi\left(\cdot\right) \right\|_2 \leq \left\| P^t\left(i, \cdot\right) - \pi\left(\cdot\right) \right\|_1 = \sum_{j=1}^{n} \left| P^t\left(i, j\right) - \pi\left(j\right) \right| \leq \frac{\varepsilon}{2},$$

где $P^t\left(i, j\right)$ – условная вероятность того, что стартуя из состояния $i \in V_G$ через $t$ шагов, марковский процесс окажется в состоянии $j \in V_G$.

В нашем случае можно явно геометрически описать множество $S$ в виде целочисленных точек множества А, попавших внутрь эллипсоида (сферы) размера $\mathrm{O}\left(\sqrt{N}\right)$ с $\pi\left(S\right) \simeq 1/2$ и с центром в точке $d^*$. На этом множестве достигается решение изопериметрической задачи в определении константы Чигера. Значение константы Чигера при этом будет пропорционально отношению площади поверхности этого эллипсоида к его объему, т.е. $h(G) \sim N^{-1/2}$. Используя это наблюдение, можно получить, что

$$T_n\left(\varepsilon, N; d\left(0\right)\right) \sim \ln\left(\pi\left(i\right)^{-1}\right) + \ln N .$$

Отсюда видно, что время выхода зависит от точки старта. Если ограничиться точками старта $i \in V_G$, для которых равномерно по $N \to \infty$ и $i, j = 1,\ldots,n$ имеет место неравенство $d_{ij}\left(0\right)/N \geq \varsigma_n > 0$, то $\ln\left(\pi\left(i\right)^{-1}\right) \sim \ln N$.

Объединяя приведенные результаты, получаем утверждение теоремы 1. ●

Эта теорема уточняет результат работы [12], явно указывая скорость сходимости и плотность концентрации. Приведенная скорость сходимости характерна для более широкого класса моделей стохастической химической кинетики, приводящих к равновесию вида неподвижной точки [5]. Плотность концентрации оценивалась на базе конструкции работы [17].

Естественно в виду теоремы 1 принимать решение задачи (2) $\left\{d_{ij}^{*}\right\}_{i=1,\, j=1}^{n,\, n}$ за равновесную конфигурацию. Обратим внимание, что предложенный выше вывод известной энтропийной модели расчета матрицы корреспонденций отличается от классического [1]. В монографии А.Дж. Вильсона [1] $\beta$ интерпретируется как множитель Лагранжа к ограничению на среднее “время в пути”: $\sum_{i,j=1}^{n} d_{ij} T_{ij} = C$.

**Замечание 1.** При этом остальные ограничения имеют такой же вид, а функционал имеет вид $F\left(d\right) = -\sum_{i,j=1}^{n} d_{ij} \ln\left(d_{ij}\right)$. Тогда, согласно экономической интерпретации двойственных множителей Л.В. Канторовича: $\beta(C) = \partial F\left(d\left(C\right)\right)/\partial C$. Из такой интерпретации иногда делают вывод о том, что $\beta$ можно понимать, как цену единицы времени в пути. Чем больше $C$, тем меньше $\beta$.

Приведенный нами вывод позволяет контролировать знак параметра $\beta > 0$ и лучше понимать его физический смысл.

**Замечание 2.** Отметим, что также как и в [1] из принципа ле Шателье–Самуэльсона следует, что с ростом $\beta$ среднее время в пути $\sum_{i,j=1}^{n} d_{ij}\left(\beta\right) T_{ij}$ будет убывать. В связи с этим обстоятельством, а также исходя из соображений размерности, вполне естественно понимать под $\beta$ величину, обратную к характерному (среднему) времени в пути [1] – физическая интерпретация. Собственно, такая интерпретация параметра $\beta$, как правило, и используется в многостадийных моделях (см., например, [18]).

В дальнейшем, нам будет удобно привести задачу (2) к следующему виду (при помощи метода множителей Лагранжа [19], теоремы фон Неймана о минимаксе [20] и перенормировки $d := d/N$):

$$\max_{\lambda^L,\lambda^W} \min_{\sum_{i,j=1}^{n} d_{ij}=1,\, d_{ij}\ge 0} \left[ \sum_{i,j=1}^{n} d_{ij}\ln d_{ij} + \beta \sum_{i=1,\, j=1}^{n,n} d_{ij}T_{ij} + \sum_{i=1}^{n} \lambda_i^L \left( \sum_{j=1}^{n} d_{ij} - l_i \right) + \sum_{i=1}^{n} \lambda_j^W \left( \mathrm{w}_j - \sum_{i=1}^{n} d_{ij} \right) \right], \quad (3)$$

где $l = L/N$, $\mathrm{w} = W/N$.

**Замечание 3.** Используя принцип Ферма [19] не сложно проверить, что решение задачи (3) можно представить в виде

$$d_{ij} = \exp\left(-\lambda_i^L\right)\exp\left(\lambda_j^w\right)\exp\left(-\beta T_{ij}\right),$$

где

$$\exp\left(\lambda_i^L\right) = \frac{1}{l_i}\sum_{j=1}^{n}\exp\left(\lambda_j^w\right)\exp\left(-\beta T_{ij}\right), \; \exp\left(-\lambda_j^w\right) = \frac{1}{w_j}\sum_{i=1}^{n}\exp\left(-\lambda_i^L\right)\exp\left(-\beta T_{ij}\right). \quad (4)$$

Отсюда можно усмотреть интерпретацию двойственных множителей как соответствующих “потенциалов притяжения/отталкивания районов” [1]. К этому мы еще вернемся в п. 3.

**Замечание 4.** Департамент транспорта г. Москвы несколько лет назад поставил следующую задачу. Сколько человек надо обзвонить и опросить на предмет того какой корреспонденции они принадлежат (где живут и где работают), чтобы восстановить матрицу корреспонденций с достаточной точностью и доверительным уровнем? Формализуем задачу. Предположим, что истинная (пронормированная) матрица корреспонденций $\left\{d_{ij}^*\right\}_{i=1,\, j=1}^{n,\, n}$ – ее и надо восстановить. В результате опросов населения получилась матрица (вектор) $r = \left\{r_{ij}\right\}_{i=1,\, j=1}^{n,\, n}$, где $r_{ij}$ – количество опрошенных респондентов, принадлежащих корреспонденции $(i, j)$, $\sum_{i,j=1}^{n} r_{ij} = N$. Задачу формализуем следующим образом: найти наименьшее $N$, чтобы

$$P_{d^*}\left(\left\|\frac{r_{ij}}{N} - d^*\right\|_2 \ge \varepsilon\right) \le \sigma.$$

Нижний индекс $d^*$ у вероятности означает, что случайный вектор $r$ имеет мультиномиальное распределение с параметром $d^*$, т.е. считаем, что опрос проводился в идеальных условиях. Из теоремы 1 не сложно усмотреть, что достаточно опросить $N = \left(4 + 8\ln\left(\sigma^{-1}\right)\right)\varepsilon^{-2}$ респондентов. Скажем, опрос ~50 000 респондентов, который и был произведен, действительно позволяет не плохо восстановить матрицу корреспонденций $d^* \approx \bar{d} \overset{def}{=} r_{ij}/N$. Однако мы привели здесь это замечание для других целей. В ряде работ (см., например [21, 22]) исходя из данных таких опросов также восстанавливают матрицу корреспонденций, но при другой параметрической гипотезе (меньшее число параметров): $d_{ij} = \exp\left(-\lambda_i^L\right)\exp\left(\lambda_j^w\right)\exp\left(-\beta T_{ij}\right)$. То есть дополнительно предполагают, что $n^2$ неизвестных параметров в действительности однозначно определяются $2n$ параметрами (иногда к ним добавляют еще один параметр $\beta$), которые у нас ранее (замечание 3) интерпретировались как множители Лагранжа. Встает вопрос: как оптимально оценить эти параметры? Ответ дает теорема Фишера (в современном не асимптотическом варианте изложение этой теоремы можно найти в [23]) об оптимальности оценок максимального правдоподобия (ОМП). Собственно, для выборки из мультиномиального распределения оценкой максимального правдоподобия как раз и будет выборочное среднее $\bar{d}$. Для описанной параметрической модели (с $2n$ параметрами) поиск такой оценки сводится разрешению системы (4) (замечание 3).

**Замечание 5.** Подобно тому, как мы рассматривали в этом пункте трудовые корреспонденции (в утренние и вечерние часы более 70% корреспонденций по Москве и области именно такие), можно рассматривать перемещения, например, из дома к местам учебы, отдыха, в магазины и т.п. (по хорошему, еще надо было учитывать перемещения типа работа–магазин–детский_сад–дом) – рассмотрение всего этого вкупе приведет также к задаче (2). Только будет больше типов корреспонденций $d$ : помимо пары районов, еще нужно будет учитывать тип корреспонденции [1]. Все это следует из того, что инвариантной мерой динамики с несколькими типами корреспонденций по-прежнему будет прямое произведение пуассоновских мер. Другое дело, когда мы рассматриваем разного типа пользователей транспортной сети, например: имеющих личный автомобиль и не имеющих личный автомобиль. Первые могут им воспользоваться равно, как и общественным транспортом, а вторые нет. И на рассматриваемых масштабах времени пользователи могут менять свой тип. То есть время в пути может для разных типов пользователей быть различным [1]. Считая, подобно тому как мы делали раньше, что желание пользователей корреспонденции $(i,j)$ сети сменить свой тип (вероятность в единицу времени) есть

$$\tilde{\lambda}\exp\Big(\underbrace{\tilde{R}\big(T_{ij}^{old}\big)}_{\text{суммарные затраты до смены типа}} - \underbrace{\tilde{R}\big(T_{ij}^{new}\big)}_{\text{суммарные затраты после смены типа}}\Big), \text{ где } \tilde{R}(T)=\beta T,$$

и учитывая в “обменах” тип пользователя (будет больше типов корреспонденций $d$, но “меняются местами работы” только пользователи одного типа), можно показать, что все это вкупе приведет также к задаче типа (2).

## 3. Вывод на основе модели равновесного распределения потоков

Предварительно напомним, следуя [11, 13] эволюционный вывод модели равновесного распределения потоков [12, 13, 24, 25].

Задан ориентированный граф $\Gamma=(V,E)$ – транспортная сеть города ($V$ – узлы сети (вершины), $E\subset V\times V$ – дуги сети (ребра графа)). В графе имеются две выделенные вершины. Одна из вершин графа является источником, другая стоком. Из источника в сток ведет много путей, которые мы будем обозначать $p\in \mathrm{P}$, $|\mathrm{P}|=m$. Обозначим через

$x_p$ – величина потока по пути $p$, $x_p\in S_m(N)=\left\{x=\{x\}_{p\in\mathrm{P}}\ge 0:\ \sum_{p\in\mathrm{P}}x_p=N\right\}$;

$y_e$ – величина потока по ребру $e\in E$: $y_e=\sum_{p\in\mathrm{P}}\delta_{ep}x_p$ ($y=\Theta x$, $\Theta=\{\delta_{ep}\}_{e\in E,p\in\mathrm{P}}$), где

$$\delta_{ep}=\begin{cases}1, & e\in p\\ 0, & e\notin p\end{cases};$$

$\tau_e(y_e)$ – удельные затраты на проезд по ребру $e$ (гладкие неубывающие функции);

$G_p(x)=\sum_{e\in E}\tau_e(y_e)\delta_{ep}$ – удельные затраты на проезд по пути $p$.

Рассмотрим, следуя [11], популяционную игру в которой имеется набор $N$ однотипных игроков (агентов). Множеством чистых стратегий каждого такого агента является $\mathrm{P}$, а выигрыш (потери со знаком минус) от использования стратегии $p\in\mathrm{P}$ определяются формулой $-G_p(x)$.

Опишем динамику поведения игроков. Пусть в момент времени $t \ge 0$ агент использует стратегию $p \in \mathrm{P}$, $\lambda_{p,q}(t)\Delta t + o(\Delta t)$ – вероятность того, что он поменяет свою стратегию на $q \in \mathrm{P}$ в промежутке времени $(t, t+\Delta t)$. Будем считать, что

$$\lambda_{p,q}(t) \equiv \lambda_{p,q} = \lambda P_q\left(\left\{G_p\left(x(t)\right)\right\}_{p\in\mathrm{P}}\right),$$

где (как и в п. 2) коэффициент $0 < \lambda = \mathrm{O}(1)$ характеризует интенсивность “перескоков” агентов, а

$$P_q\left(\left\{G_p\left(x(t)\right)\right\}_{p\in\mathrm{P}}\right) = P\left(q = \arg\max_{p\in\mathrm{P}}\left\{-G_p\left(x(t)\right) + \xi_p\right\}\right).$$

Если $\xi_p \equiv 0$, то получаем динамику наилучших ответов [11], если $\xi_p$ – независимые одинаково распределенные случайные величины с распределением Гумбеля [26]: $P\left(\xi_p < \xi\right) = \exp\left\{-e^{-\xi/\omega - E}\right\}$, где $\omega \in (0, \omega_0]$ ($\omega_0 = \mathrm{O}(1)$), $E \approx 0.5772$ – константа Эйлера, а $\mathrm{Var}\left[\xi_p\right] = \omega^2\pi^2/6$, то получаем логит динамику [11]

$$P_q\left(\left\{G_p\left(x(t)\right)\right\}_{p\in\mathrm{P}}\right) = \frac{\exp\left(-G_q\left(x(t)\right)/\omega\right)}{\sum_{p\in\mathrm{P}}\exp\left(-G_p\left(x(t)\right)/\omega\right)},$$

вырождающуюся в динамику наилучших ответов при $\omega \to 0+$. Далее мы будем считать, что задана логит динамика. Объясняется такая динамика совершенно естественно. Каждый агент имеет какую-то картину текущего состояния загрузки системы, но либо он не имеет возможности наблюдать ее точно, либо он старается как-то спрогнозировать возможные изменения загрузок (а соответственно и затраты на путях) в будущем (либо и то и другое).

Приведем соответствующий аналог теоремы 1.

**Теорема 2.** *Для любого* $x(0) \in S_m(N)$ *существует такая константа* $c_m\left(x(0)\right) > 0$, *что для всех* $\sigma \in (0, 0.5)$, $t \ge c_m\left(x(0)\right)\ln N$ *имеет место неравенство*

$$P\left(\left\|\frac{x(t)}{N} - x^*\right\|_2 \ge \frac{2\sqrt{2} + 4\sqrt{\ln\left(\sigma^{-1}\right)}}{\sqrt{N}}\right) \le \sigma,\ x(t) \in S_m(N),$$

*где (стохастическое равновесие Нэша–Вардропа)*

$$x^* = \arg\min_{x\in S_m(1)} \Psi\left(y(x)\right) + \omega\sum_{p\in\mathrm{P}} x_p \ln\left(x_p\right), \tag{5}$$

$$\Psi\left(y(x)\right) = \sum_{e\in E}\int_0^{y_e(x)} \tilde{\tau}_e(z)\,dz,\ \tilde{\tau}_e(z) = \tau_e(zN).$$

**Схема доказательства.** Стационарная мера описанного марковского процесса имеет (с точностью до нормирующего множителя) вид (теорема 11.5.12 [11])

$$\frac{N!}{x_1!\cdot\ldots\cdot x_m!}\exp\left(-\frac{\Psi\left(y\left(x\right)\right)}{\omega}\right),\ x\in S_m\left(N\right).$$

Для оценок плотности концентрации достаточно заметить, что $\Psi\left(y\left(x\right)\right)$ – выпуклая функция, как композиция линейной и выпуклой функции. Поэтому оценки плотности концентрации здесь не могут быть хуже оценок в теореме 1. Рассуждения для оценки времени выхода проводятся аналогично теореме 1. ●

**Следствие (Proposition 1 [27]).** *Если $\omega\to 0+$, то решение задачи (5) сходится к*

$$x^* = \arg\min_{x\in S_m(1):\ \Theta x = y^*}\sum_{p\in\mathrm{P}} x_p\ln\left(x_p\right),\tag{6}$$

*где* $y^* = \arg\min_{y=\Theta x,\ x\in S_m(1)}\Psi\left(y\right)$.

Это следствие решает проблему обоснования гипотезы Бар-Гира [28]. Напомним вкратце в чем состоит эта гипотеза. Известно, см., например, [25], что поиск равновесного распределения потоков по путям (равновесия Нэша–Вардропа) в модели с одним источником и стоком сводится к задаче оптимизации (5) с $\omega=0$. Хотя функционал этой задачи выпуклый, но он не всегда строго выпуклый, даже в случае, когда все функции $\tau_e\left(y_e\right)$ – строго возрастающие (тогда можно лишь говорить о единственности равновесного распределения потоков по ребрам $y^*$). Гипотеза Бар-Гира говорит, что "в жизни" с большой вероятностью реализуется то равновесие из множества равновесий, которое находится из решения задачи (6).

**Замечание 6.** Распределение потоков по путям $x$ называется равновесием (Нэша–Вардропа) в рассматриваемой популяционной игре $\left\langle\left\{x_p\right\}_{p\in\mathrm{P}},\left\{G_p\left(x\right)\right\}_{p\in\mathrm{P}}\right\rangle$, если из $x_p>0$, $p\in\mathrm{P}$ следует $G_p\left(x\right)=\min_{q\in\mathrm{P}}G_q\left(x\right)$.

Или, что то же самое:

для любых $p\in\mathrm{P}$ выполняется $x_p\cdot\left(G_p\left(x\right)-\min_{q\in\mathrm{P}}G_q\left(x\right)\right)=0.$

**Замечание 7.** Если при $\omega=0$ рассмотреть предельный случай

$$\tau_e\left(y_e\right):=\tau_e^{\mu}\left(y_e\right)\xrightarrow[\mu\to 0+]{}T_e,$$

то поиск равновесия Нэша–Вардропа просто сводится к поиску социального оптимума, что приводит в данном случае к решению транспортной задачи линейного программирования [29]. Если делать предельный переход с учетом ограничений на пропускные способности ребер

$$\tau_e\left(y_e\right):=\tau_e^{\mu}\left(f_e\right)\xrightarrow[\mu\to 0+]{}\begin{cases}T_e, & 0\le y_e<\bar{y}_e\\ \left[T_e,\infty\right), & y_e=\bar{y}_e\end{cases},$$

то получится более сложная задача, которая описывает равновесие в (стохастической, если $\omega>0$) модели стабильной динамики [13, 30].

Из рассмотренных в замечании 7 случаев выпала ситуация $\omega > 0$ $\tau_e(y_e) := \tau_e^{\mu}(y_e) \xrightarrow[\mu \to 0+]{} T_e$. Ее мы сейчас отдельно и рассмотрим на примере другого эволюционного способа обоснования энтропийной модели расчета матрицы корреспонденций. По ходу обсуждения этого способа с нашими коллегами (прежде всего, Ю.Е. Нестеровым и Ю.В. Дорном) у предложенного подхода появилось название: “облачная модель”.

Предположим, что все вершины, отвечающие источникам корреспонденций, соединены ребрами с одной вспомогательной вершиной (облако № 1). Аналогично, все вершины, отвечающие стокам корреспонденций, соединены ребрами с другой вспомогательной вершиной (облако № 2). Припишем всем новым ребрам постоянные веса. И проинтерпретируем веса ребер, отвечающих источникам $\lambda_i^L$, например, как средние затраты на проживание (в единицу времени, скажем, в день) в этом источнике (районе), а веса ребер, отвечающих стокам $\lambda_j^W$, как уровень средней заработной платы (в единицу времени) в этом стоке (районе), если изучаем трудовые корреспонденции. Будем следить за системой в медленном времени, то есть будем считать, что равновесное распределение потоков по путям стационарно. Поскольку речь идет о равновесном распределении потоков, то нет необходимости говорить о затратах на путях или ребрах детализированного транспортного графа, достаточно говорить только затратах (в единицу времени), отвечающих той или иной корреспонденции. Таким образом, у нас есть взвешенный транспортный граф с одним источником (облако 1) и одним стоком (облако 2). Все вершины этого графа, кроме двух вспомогательных (облаков), соответствуют районам в модели расчета матрицы корреспонденций. Все ребра этого графа имеют стационарные (не меняющиеся и не зависящие от текущих корреспонденций) веса $\left\{T_{ij}; \lambda_i^L; -\lambda_j^W\right\}$. Если рассмотреть естественную в данном контексте логит динамику ($x \equiv d$) или ее обобщение описанное выше с $\omega = 1/\beta$ (здесь полезно напомнить, что согласно замечанию 2 $\beta$ обратно пропорционально средним затратам, а $\omega$ имеет как раз физическую размерность затрат), то поиск равновесия рассматриваемой макросистемы согласно теореме 1 приводит (в прошкалированных переменных) к задаче, сильно похожей на задачу (3) из п. 2

$$\min_{\sum_{i,j=1}^{n} d_{ij}=1,\, d_{ij}\ge 0} \left[ \sum_{i,j=1}^{n} d_{ij} \ln d_{ij} + \beta \sum_{i=1,\, j=1}^{n,n} d_{ij} T_{ij} + \beta \sum_{i=1}^{n} \left( \lambda_i^L \sum_{j=1}^{n} d_{ij} \right) - \beta \sum_{i=1}^{n} \left( \lambda_j^W \sum_{i=1}^{n} d_{ij} \right) \right].$$

Разница состоит в том, что здесь мы не оптимизируем по $2n$ двойственным множителям $\lambda^L$, $\lambda^W$ (множителям Лагранжа). Более того, мы их и не интерпретируем здесь как двойственные множители, поскольку мы их ввели на этапе взвешивания ребер графа. Тем не менее, значения этих переменных, как правило, не откуда брать. Тем более что приведенная выше (наивная) интерпретация вряд ли может всерьез рассматриваться, как способ определения этих параметров исходя из данных статистики. Более правильно понимать $\lambda_i^L$, $\lambda_j^W$ как потенциалы притяжения/отталкивания районов (см. также замечание 3), включающиеся в себя плату за жилье и зарплату, но включающие также и многое другое, что сложно описать количественно. И здесь как раз помогает информация об источниках и стоках, содержащаяся в $2n$ уравнениях задающих множество А. Таким образом, мы приходим ровно к той же самой задаче (3) с той лишь разницей, что мы получили дополнительную интерпретацию двойственных множителей в задаче (3). При этом двойственные множители в задаче (3) равны (с точностью до мультипликативного фактора $\beta$) введенным здесь потенциалам притяжения районов.

Нам представляется очень плодотворной и перспективной идея перенесения имеющейся информации об исследуемой системе из обременительных законов сохранения динамики, описывающей эволюцию этой системы, в саму динамику путем введения дополнительных естественно интерпретируемых параметров. При таком подходе становится возможным, например, учитывать в моделях и рост транспортной сети. Другими словами, при таком подходе, например, можно естественным образом рассматривать также и ситуацию, когда число пользователей транспортной сетью меняется со временем (медленно).

Исследование выполнено в рамках федеральной целевой программы «Исследования и разработки по приоритетным направлениям развития научно-технологического комплекса России на 2014 – 2020 годы», Соглашение № 14.604.21.0052 от 30.06.2014 г. с Минобрнаукой. Уникальный идентификатор проекта RFMEFI60414X0052.

Исследование было поддержано грантом РФФИ 15-31-70001 мол_а_мос.

## СПИСОК ЦИТИРОВАННОЙ ЛИТЕРАТУРЫ

1. *Вильсон А.Дж.* Энтропийные методы моделирования сложных систем. М.: Наука, 1978.

   *A.G. Wilson* Entropy in urban and regional modeling. Routledge, 2011.

2. *Ethier N.S., Kurtz T.G.* Markov processes. Wiley Series in Probability and Mathematical Statistics: Probability and Mathematical Statistics. John Wiley & Sons Inc., New York, 1986.

3. *Малышев В.А., Пирогов С.А.* Обратимость и необратимость в стохастической химической кинетике // Успехи мат. наук. 2008. Т. 63. вып. 1(379). С. 4–36.

   *V. A. Malyshev, S. A. Pirogov*, "Reversibility and irreversibility in stochastic chemical kinetics", *Uspekhi Mat. Nauk*, 63:1(379) (2008), 3–36

4. *Батищева Я.Г., Веденяпин В.В.* II-й закон термодинамики для химической кинетики // Матем. Моделирование. 2005. Т. 17:8. С. 106–110.

   *Ya. G. Batishcheva, V. V. Vedenyapin*, "The 2-nd low of thermodynamics for chemical kinetics", *Matem. Mod.*, 17:8 (2005), 106–110

5. *Гасников А.В., Гасникова Е.В.* Об энтропийно-подобных функционалах, возникающих в стохастической химической кинетике при концентрации инвариантной меры и в качестве функций Ляпунова динамики квазисредних // Математические заметки. 2013. Т. 94. № 6. С. 816–824. arXiv:1410.3126

   *A. V. Gasnikov, E. V. Gasnikova*, "On Entropy-Type Functionals Arising in Stochastic Chemical Kinetics Related to the Concentration of the Invariant Measure and Playing the Role of Lyapunov Functions in the Dynamics of Quasiaverages", *Mat. Zametki*, 94:6 (2013), 819–827

6. *Boucheron S., Lugoshi G., Massart P.* Concentration inequalities: A nonasymptotic theory of independence. Oxford University Press, 2013.

7. *Санов И.Н.* О вероятности больших отклонений случайных величин // Матем. сб. 1957. Т. 42(84):1. С. 11–44.

   *I. N. Sanov*, "On the probability of large deviations of random magnitudes", *Mat. Sb. (N.S.)*, 42(84):1 (1957), 11–44

8. *Fang S.-C., Rajasekera J.R., Tsao H.-S.J.* Entropy optimization and mathematical programming. Kluwer's International Series, 1997.

9. *Гардинер К.В.* Стохастические методы в естественных науках. М.: Мир, 1986.

   *Gardiner C.* Stochastic methods. A Handbook for the Natural and Social Sciences. Springer, 2009.

10. *Вайдлих В.* Социодинамика: системный подход к математическому моделированию в социальных науках. М.: УРСС, 2010.

    *Weidlich W.* Sociodynamics: a System Approach to Mathematical Modellig in the Social Sciences. Amsterdam: Harwood Academic Publishers, 2000.

11. *Sandholm W.* Population games and Evolutionary dynamics. Economic Learning and Social Evolution. MIT Press; Cambridge, 2010.

12. *Гасников А.В., Кленов С.Л., Нурминский Е.А., Холодов Я.А., Шамрай Н.Б.* Введение в математическое моделирование транспортных потоков. Под ред. А.В. Гасникова с приложениями М.Л. Бланка, К.В. Воронцова и Ю.В. Чеховича, Е.В. Гасниковой, А.А. Замятина и В.А. Малышева, А.В. Колесникова, Ю.Е. Нестерова и С.В. Шпирко, А.М. Райгородского, с предисловием руководителя департамента транспорта г. Москвы М.С. Ликсутова. М.: МЦНМО, 2013. 427 стр., 2-е изд.

Introduction to the mathematical modelling of traffic flows. Eds. A.V. Gasnikov. Moscow: MCCME, 2013. [in Russian]

13. *Гасников А.В., Дорн Ю.В., Нестеров Ю.Е., Шпирко С.В.* О трехстадийной версии модели стационарной динамики транспортных потоков // Математическое моделирование. Т. 26. № 6. 2014. С. 34–70. arXiv:1405.7630

*A. Gasnikov, Yu. Dorn, Yu. Nesterov, S. Shpirko*, "On the three-stage version of stable dynamic model", *Matem. Mod.*,26:6 (2014), 34–70

14. *Гасников А.В., Гасникова Е.В., Нестеров Ю.Е., Чернов А.В.* Об эффективных численных методах решения задач энтропийно-линейного программирования в случае разреженных матриц // ЖВМиМФ. 2015. (принята к печати) arXiv:1410.7719

*A.V. Gasikov, E.V. Gasnikova, Yu.E. Nesterov, A.V. Chernov*, "Entropy-linear programming" Comp. Math. and Math. Phys. 2015. (in print)

15. *Боровков А.А.* Эргодичность и устойчивость случайных процессов. М.: УРСС, 1999.

*Borovkov A.*A. Ergodicity and Stability of Stochastic Processes. Wiley series in Probability and Statistics, 1998.

16. *Levin D.A., Peres Y., Wilmer E.L.* Markov chain and mixing times. AMS, 2009. http://pages.uoregon.edu/dlevin/MARKOV/markovmixing.pdf

17. *Гасников А.В., Дмитриев Д.Ю.* Об эффективных рандомизированных алгоритмах поиска вектора PageRank // ЖВМиМФ. 2015. Т. 54. № 3. С. 355–371. arXiv:1410.3120

*A.V. Gasnikov, D.Yu. Dmitriev*, "Efficient randomized algorithms for PageRank problem" Comp. Math. and Math. Phys. 2015. V. 54. no. 3. P. 355–371.

18. *Ortúzar J.D., Willumsen L.G.* Modelling transport. JohnWilley & Sons, 2011.

19. *Магарил-Ильяев Г.Г., Тихомиров В.М.* Выпуклый анализ и его приложения. М.: УРСС, 2011.

*G.G. Magaril-Il'yaev, V.M. Tikhomirov* Convex Analysis: Theory and Applications (Translation of Mathematical Monographs). AMS, V. 222. 2003.

20. *Sion M.* On general minimax theorem // Pac. J. Math. 1958. V. 8. P. 171–176.

21. *Yun S., Sen A.* Computation of maximum likelihood estimates of gravity model parameters // Journal of Regional Science. 1994. V. 34. № 2. P. 199–216.

22. *Sen A.* Maximum likelihood estimation of gravity model parameters // Journal of Regional Science. 1986. V.26. № 3. P. 461–474.

23. *Spokoiny V.* Parametric estimation. Finite sample theory // The Annals of Statistics. 2012. V. 40. № 6. P. 2877–2909. arXiv:1111.3029v5

24. *Sheffi Y.* Urban transportation networks: Equilibrium analysis with mathematical programming methods. N.J.: Prentice–Hall Inc., Englewood Cliffs, 1985.

25. *Patriksson M.* The traffic assignment problem. Models and methods. Utrecht, Netherlands: VSP, 1994.

26. *Andersen S.P., de Palma A., Thisse J.-F.* Discrete choice theory of product differentiation. MIT Press; Cambridge, 1992.

27. *Cuturi M., Peyré G., Rolet A.* A smoothed approach for variational Wasserstein problem // e-print, 2015. arXiv:1503.02533

28. *Bar-Gera H.* Origin-based algorithms for transportation network modeling. Univ. of Illinois at Chicago, 1999.

29. *Ahuja R.K., Magnati T.L., Orlin J.B.* Network flows: Theory, algorithms and applications. Prentice Hall, 1993.

30. *Nesterov Y., de Palma A.* Stationary Dynamic Solutions in Congested Transportation Networks: Summary and Perspectives // Networks Spatial Econ. 2003. № 3(3). P. 371–395.